\documentclass[11pt, a4paper,reqno]{amsart}

\usepackage{amscd,amsmath}
\usepackage{amssymb,amsmath,latexsym,color,enumerate,tikz,dsfont}
\usepackage[all]{xypic}       % Package XYpic of commutative diagram macros
\usepackage{mathabx}
\usepackage{mathrsfs}
\usepackage{stmaryrd}
\usepackage{setspace}
\usepackage[makeroom]{cancel}

\usepackage{tikz}
\usetikzlibrary{matrix}

\usepackage{floatrow}

\def\2{\text{\bf 2}}

\def\Frm{{\mathsf{Frm}}}
\def\Top{{\mathsf{Top}}}

\newcommand{\newuparrow}{{{\rlap{$\ $}\hbox{$\uparrow$}}}}%
\newcommand{\twoheaddownarrow}{{\rlap{\rlap{$\ $}\raise .25ex\hbox{$\downarrow$}}\raise-.25ex\hbox{$\downarrow$}}}
\newcommand{\twoheaduparrow}{{\rlap{\rlap{$\ $}\raise .25ex\hbox{$\uparrow$}}\raise-.25ex\hbox{$\uparrow$}}}

 \newcommand{\tbigcup}{\mathop{\textstyle \bigcup }}%%
 \newcommand{\tbigvee
}{\mathop{\textstyle \bigvee }}%%
 \newcommand{\tbigwedge
}{\mathop{\textstyle \bigwedge }}%%

\newcommand{\dirvee}{\mathop{\setminus\kern-.85ex\nnearrow}}

\DeclareMathAccent{\widetriangle}{\mathord}{largesymbols}{"E6}

\newtheorem{theorem}{Theorem}[section]
\newtheorem{proposition}[theorem]{Proposition}
\newtheorem{lemma}[theorem]{Lemma}
\newtheorem{corollary}[theorem]{Corollary}

\theoremstyle{definition}
\newtheorem{definition}[theorem]{Definition}
\newtheorem{example}[theorem]{Example}

\theoremstyle{remark}

\newtheorem{remark}[theorem]{Remark}

\newcommand{\msum}{\mathrel{\mathcal{S}}}

\begin{document}

\title[The Boolean completion of $\mathsf{S}(L)$]{The Boolean completion of $\mathsf{S}(L)$:\\  a refinement of localic mereology}

%% First author
\author[M.\thinspace A. Moshier]{M. Andrew Moshier}
\address[M. Andrew Moshier]{CECAT\\ Department of Mathematics \& Computer Science\\ Chapman University\\ Orange CA 92866\\ U.S.A.}
\email{moshier@chapman.edu}
%
%%% Second author
\author[I. Mozo~Carollo]{Imanol Mozo~Carollo}
\address[Imanol Mozo~Carollo]{Department of Applied Economics\\ University of the Basque Country UPV/EHU\\ Apdo. 644\\ 48080 Bilbao\\ Spain}
\email{imanol.mozo@ehu.eus}
%
%% Third author
%\author[J. Walters-Wayland]{Joanne Walters-Wayland}
%\address[Joanne Walters-Wayland]{CECAT\\ Department of Mathematics \& Computer Science\\ Chapman University\\ Orange CA 92866\\ U.S.A.}
%\email{joanne@waylands.com}
%
%
%
%\author[J. Picado]{Jorge~Picado}
%\address{CMUC\\Department of Mathematics\\University of Coimbra\\3001-501 Coimbra\\Portugal}
%\email{picado@mat.uc.pt}
%\urladdr{http://www.mat.uc.pt/~picado}

\date{\today}



\thanks{The second named author acknowledges support from the Basque Government (grant IT1913-26).\\
Accepted in Mathematica Slovaca}

%\setstretch{1.3}

\maketitle
%\tableofcontents

\begin{abstract}
We construct and study the Boolean completion of the coframe $\mathsf{S}(L)$ of sublocales of a frame $L$. Realizing $\mathsf{S}(L)$ as a basis for an Alexandroff topology on $L$ and passing to the regular elements of the associated frame, we obtain a complete Boolean algebra $\mathsf{P}(L)$ whose elements are certain subsets of $L$ generalizing sublocale sets; we call them the \emph{parts} of $L$. We show that $\mathsf{S}(L)$ embeds into $\mathsf{P}(L)$ as the left adjoint of a Galois connection and that the canonical embedding of $L$ into $\mathsf{S}(L)$ extends to $\mathsf{P}(L)$. From a mereological standpoint, this completion is motivated by universalism: after proving that mereological sums and exact joins coincide for separative posets, we identify $\mathsf{P}(L)$ as the sum-completion (equivalently, the  exact-join-completion) of $\mathsf{S}(L)$, thereby supplying the mereological sums that $\mathsf{S}(L)$ lacks. Finally, we single out several distinguished parts of a frame: the pointless part, which extends Ball's construction beyond the $T_1$ setting, together with the part of merely points and the topological glue, the former  of the two having no counterpart in the classical theory of sublocales.
\end{abstract}

\section*{Introduction}

Sublocale sets are a particularly lucid way of representing ``subspaces'' in pointfree topology for several reasons: their set-theoretic inclusion represents their geometrical order,  lattice operations are especially simple to compute (meets are intersections and joins are the closures under subset infima of the unions), there is even a nice parallel between the relation of points in a space and the corresponding one-point subspaces on one side and points in a frame and one-point sublocales on the other \cite{PP2008}\dots. In precise terms, a 
\emph{sublocale set} (\emph{sublocale} for short) of a frame $L$ is a subset $S$ of $L$  satisfying
\begin{itemize}
	\item[(S1)] for every $A \subseteq S$, $\bigwedge A$ is in $S$, and
	\item[(S2)] for
	every $s \in S$ and every $x \in L$, $x \to s$ is in $S$.
\end{itemize}

In the classical setting, a space is determined by its points, and so are its subspaces. Consequently, the poset of subspaces of a space is just the powerset of its set of points, hence a Boolean algebra. This picture changes in the pointfree setting. For a frame $L$, the collection $\mathsf{S}(L)$ of sublocale sets is far richer than a powerset. Ordered by inclusion, $\mathsf{S}(L)$ forms a coframe, which is in general far from being Boolean \cite{PP12}. This non-Boolean behaviour leads to distinctive phenomena, such as Isbell’s Density Theorem \cite{I72}, which states that every locale has a least dense sublocale, or the avoidance of the Banach-Tarski paradox \cite{S2012}.

The main goal of this paper is to construct and study the Boolean completion of $\mathsf{S}(L)$ in terms of certain subsets of $L$ that naturally generalize sublocale sets. The motivation for this work is twofold. One part is purely algebraic  and requires no further justification beyond its intrinsic mathematical interest. The other is more conceptual and we explain it below.

The theory of sublocales provides a framework for analysing the \emph{mereological} structure of spatial or geometrical objects. Mereology studies parthood relations---relations between part and whole, and between parts within a whole---in a very general sense. Its formal development goes back to the work of Stanis\l aw Le\'sniewski \cite{L16, L27}; for modern systematic treatments see Simons \cite{S87} and Cotnoir-Varzi \cite{CV21}. From a mereological perspective, $\mathsf{S}(L)$ fails to satisfy a classical principle known as \emph{universalism}, which informally states that any nonempty collection of parts of a whole can be composed into a part. At first sight, at least for a mathematician, this may seem unproblematic: $\mathsf{S}(L)$ is a complete lattice, so any family of sublocales has a supremum. However, in the mereological tradition, composition is not identified with taking suprema, although this view has become more common in recent work.

A mereological composition is meant to be a part of the whole that is \emph{entirely composed} of the given parts. More precisely, it must contain all the parts in the collection, and any further part of the composition must overlap one of them. We refer to the formal counterpart of this idea as a \emph{sum}. In \cite{Mozo2024}, we studied how to extend mereological models so as to satisfy universalism, and showed that, assuming a standard principle---strong supplementation---there is a unique such extension, namely the Boolean completion of the original model \cite{Jech2003}. Since $\mathsf{S}(L)$, viewed as a mereological model (to be more precise, $\mathsf{S}(L)$ with the bottom element removed, since the original approaches to formal mereology avoid the notion of an empty part owing to their nominalist commitments), satisfies strong supplementation (equivalently, it is a separative poset) it is natural to ask what its Boolean completion looks like.

After exploring the relation between mereological sums and exact joins and showing in Section \ref{sums.exactjoins} that for separative posets both notions coincide, we construct in Section \ref{section.BooleanCompletion} the Boolean completion of $\mathsf{S}(L)$ as follows: take the topology on $L$ generated by $\mathsf{S}(L)$, drop the empty-set, and take the complete Boolean algebra of regular elements of the resulting frame. Sublocale sets turn out to be elements of this complete Boolean algebra and this allows us to consider the elements of the Boolean completion as ``generalized sublocale sets''. In Section \ref{section.Parts} we refer to the elements of this Boolean completion as \emph{parts} of the original frame and denote the complete Boolean algebra they form by $\mathsf{P}(L)$. We also show that there is a Galois connection between $\mathsf{S}(L)$ and $\mathsf{P}(L)$ where the inclusion of $\mathsf{S}(L)$ into $\mathsf{P}(L)$ is the left adjoint. Further, the frame embedding of $L$ into $\mathsf{S}(L)$ which maps each element $a$ of $L$ to the corresponding open sublocale $\mathfrak{o}(a)$ extends naturally to a frame embedding into $\mathsf{P}(L)$.

When \(L\) is spatial, that is, when $L \simeq \mathcal{O}X$ for some topological space $X$, these parts admit a natural interpretation as generalized subspaces of $X$. This is, in fact, a generalization of a generalization, since sublocales themselves are commonly described as \emph{generalized subspaces}. Each subspace $A \subseteq X$ induces a sublocale $\widetilde{A}$ of $L$. If $X$ is a $T_D$-space, this assignment provides a faithful representation of subspaces: for all $A,B \subseteq X$,
\[
A \subseteq B \iff \widetilde{A} \subseteq \widetilde{B}.
\]
In particular, under the $T_D$-axiom, the map
\[
\begin{aligned}
2^X& \longrightarrow \mathsf{S}(L)\\
 A &\longmapsto \widetilde{A}
\end{aligned}
\]
is injective. Moreover, this map preserves finite joins, a fact that holds even without assuming the $T_D$ axiom. As every sublocale can be expressed as the infimum of a finite join of an open and a closed sublocale, and both open and closed sublocales are induced by open and closed subspaces, it follows that all sublocales arise as infima of subspace-induced sublocales. From this perspective, the passage from subspaces of $X$ to sublocales amounts to enlarging the notion of subspace by redefining what intersections are. As shown in \cite{Mozo2024}, extending $\mathsf{S}(L)$ to its Boolean completion $\mathsf{P}(L)$  amounts to adding the missing mereological sums, which turn out to be precisely the \emph{exact joins}. Accordingly, $\mathsf{P}(L)$ can be understood as the exact-join-completion of $\mathsf{S}(L)$. So in this sense, in this second step, we further enlarge the notion of subspace by redefining what unions are.

Section \ref{section.interesting} illustrates the usefulness of our approach by identifying several distinguished parts of a frame. We recover the largest pointless part studied by Ball in \cite{B24}, extending the notion beyond $T_1$-frames, and introduce a further canonical piece: the \emph{part of merely points}. Unlike the pointless part, the  part of merely points has no analogue in the classical theory of sublocales. The merely points part provides a canonical way of separating the points into an independent component of the spatial reflection, while the pointless part of the largest spatial sublocale measures the residual interaction between the spatial and pointless aspects of the frame.

\section{Background}

\subsection{Notation}
Given a poset $(X,\leq)$ and a subset $A\subseteq X$, we denote by $\tbigvee A$ its
join (supremum) whenever it exists, abbreviating $\tbigvee\{a,b\}$ as $a\vee b$;
dually, meets (infima) are written $\tbigwedge A$ and $a\wedge b$. The upset
$\{x\mid x\geq a\}$ is denoted by ${\uparrow}a$. When they exist, the least element
of the poset (namely $\tbigvee\varnothing$) is written $0$ and the greatest one
(namely $\tbigwedge\varnothing$) is written $1$.

Let $L$ be a poset and $a,b \in L$. Following \cite{B84}, we will use the following notation:
\[
a\downarrow b=\{x\mid x\wedge b\leq a\}.
\]

If $P$ is a poset  with a bottom element $0$, we will denote $P^+=P\setminus\{0\}$. If $P$ has no bottom element, $P^+=P$.

\subsection{Posets and mereology}
If $(P,\leq)$ is a poset, we will say that $a,b\in P$ \emph{overlap} and denote it by $a\circ b$ if there exists $c\in P$ such that $c\leq a,b$. If $a,b\in P$ do not overlap, we will say they are \emph{disjoint} and denote it by $a\wr b$.

A poset $P$ is said to be \emph{separative}  if for all $a,b\in P$, if $a\not\leq b$ then there exists $c\in P$ such that $c\leq a$ and $c\wr b$.  Note that the only separative poset with a bottom element $0$ is the poset with  the singleton  $\{0\}$ as its underlying set.

The condition just introduced is well known in mereology under a different name: \emph{strong supplementation} is precisely the principle asserting that if $x$ is not
a part of $y$, then $x$ has a part disjoint from $y$ \cite{CV21}. Thus if we interpret the order relation of a poset $P$ as a parthood relation, then $P$ is separative exactly when the parthood relation satisfies strong supplementation.

\subsection{Boolean completion}
Every separative poset $(P,\leq)$ embeds into a complete Boolean algebra
$\mathbb{B}(P)$, its \emph{Boolean completion}, which is unique up to isomorphism and
in which the image of $P$ is \emph{dense}: if we identify $P$ with its image in $\mathbb{B}(P)$, for each  $x\in\mathbb{B}(P)^+$ 
there exists some $p\in P$ such that $x\leq p$ \cite{Jech2003}.

\subsection{Exact joins}
Recall from \cite{B84} and \cite{BPP13} that an element $b$ is said to be the \emph{exact join} of a subset $A$ of a lattice $L$ if
\begin{itemize}
\item $b$ is an upper bound of $A$, and 
\item if $A\subseteq c\downarrow d$ implies $b\in c\downarrow d$ for all $c<d$.
\end{itemize}

\begin{proposition}[{\cite[Proposition 2.4.1]{BPP13}}]
Exact joins are joins.
\end{proposition}

\begin{theorem}[See {\cite[Theorem 2.4.2]{BPP13} and \cite[Lemma 1.8]{B84}}]\label{exactdist}
	
Let $b$ be an element and $A$ a subset of a lattice $L$. Then the following are equivalent:
\begin{enumerate}[\rm(1)]
	\item $b$ is the exact join of $A$;
	\item $\tbigvee_{a\in A}(a\wedge x)=b\wedge x$ for all $x$.
	\end{enumerate}
\end{theorem}

\subsection{Sums} Let $P$ be a poset, $b\in P$ and $A\subseteq P$. We say that $b$ is a \emph{sum} of $A$ if 
\begin{itemize}
	\item $b$ is an upper bound of $A$, and
	\item for all $x\leq b$, there exists $a\in A$ such that $a\circ x$.
\end{itemize}
We will denote by $a\msum B$ the fact that $a\in P$ is a sum of a subset $B\subseteq P$. Note that if $a\msum B$ then $B\neq\varnothing$.

\subsection{Frames}

Our reference for the general theory of frames is Picado--Pultr \cite{PP12};
we only recall here the notions strictly needed below.

A \emph{frame} (equivalently, a \emph{locale}) is a complete lattice $L$ in which
finite meets distribute over arbitrary joins, that is,
\[
a\wedge \tbigvee B=\tbigvee \{a\wedge b\mid b\in B\}
\]
holds for every $a\in L$ and every $B\subseteq L$, and a \emph{frame homomorphism} is a map $h\colon L\to M$
preserving finite meets (the top $1$ included) and arbitrary joins (the bottom $0$
included); frames together with these maps form the category $\Frm$.

Each frame $L$  is a complete
Heyting algebra: its Heyting arrow $\rightarrow$ is characterized by
\[
a\wedge b\le c \Leftrightarrow a\le b\rightarrow c.
\]
The operation $\rightarrow $ and some of its basic properties (e.g. $a\rightarrow b=1$ iff $a\leq b$, $1\rightarrow a=a$, $a\rightarrow (b\rightarrow c)=(a\wedge b)\rightarrow c)$, and $\tbigvee A\rightarrow b=\tbigwedge_{a\in A}(a\rightarrow b)$) will be used in the sequel.

 The
\emph{pseudocomplement} of $a\in L$ is
\[
a^{\ast}=a\rightarrow 0=\tbigvee \{b\in L\mid a\wedge b=0\},
\]
and $a$ is said to be \emph{complemented} when $a\vee a^\ast=1$, and \emph{regular}
when $a^{\ast\ast}=a$ (equivalently, when $a$ is a pseudocomplement, i.e. $a=b^\ast$
for some $b$).

The motivating example is the frame $\Omega X$ of open sets of a topological space
$X$. Since taking inverse images is functorial, this assignment yields a
contravariant functor $\Omega\colon\Top\to\Frm$, where $\Top$ is the category of
topological spaces and continuous maps. In the reverse direction lies the
\emph{spectrum functor} $\Sigma\colon\Frm\to\Top$: it sends a frame $L$ to its
spectrum $\Sigma L$, whose points are the homomorphisms $\xi\colon L\to\{0,1\}$ and
whose open sets are the $\Sigma_a=\{\xi\in\Sigma L\mid \xi(a)=1\}$ with $a\in L$, and
it sends a homomorphism $h\colon L\to M$ to the continuous map
$\Sigma h\colon\Sigma M\to\Sigma L$, $\Sigma h(\xi)=\xi h$. This functor is right
adjoint to $\Omega$, the unit and counit being
$\eta_L\colon L\to\Omega\Sigma L$, $\eta_L(a)=\Sigma_a$, and
$\epsilon_X\colon X\to\Sigma\Omega X$, $\epsilon_X(x)=\hat{x}$ with
$\hat{x}(U)=1$ precisely when $x\in U$; the former is the \emph{spatial reflection}
of $L$. A frame is \emph{spatial} if it arises, up to isomorphism, as the open-set
frame of some space.

The particular notions we will need are the following:
a frame $L$ is
 \begin{enumerate}[-]
 %\item \emph{regular} if $a = \tbigvee\{b\in L\mid b \prec a\}$ for every $a\in L$, where $b\prec a$ ($b$ is \emph{rather below} $a$) means that $b^\ast\vee a=1$;

\item \emph{Boolean} if it is a Boolean algebra, that is, if every element is complemented;

\item \emph{subfit} if for each $a,b\in L$ such that $a\not\leq b$ there exists $c\in L$ such that $a\vee c=1\neq b\vee c$; and
 
% \item \emph{completely regular} if $a = \tbigvee\{b\in L\mid b{\prec\hskip-2pt \prec\hskip2pt}a\}$ for each $a\in L$, where $b{\prec\hskip-2pt \prec\hskip2pt}a$ ($b$ is \emph{completely below} $a$) means that there is $\{c_{r}\mid r\in \mathbb{Q}\cap[0,1]\}\subseteq L$ such that $b\le c_{0}$, $c_{1}\le a$ and $c_r\prec c_s$ (i.e. $c_{r}^\ast\vee c_{s}=1$) whenever $r<s$;
  
 % \item\emph{compact} if for each $A\subseteq L$ such that $\tbigvee A=1$ there exists a finite $F\subseteq A$ such that $\tbigvee F=1$;
 
% \item\emph{continuous} if $a = \tbigvee\{b\in L\mid b\ll a\}$ for every $a\in L$, where $b\ll a$ ($b$ is \emph{way below} $a$) means that ${a}\leq \tbigvee A$ for some $A\subseteq L$ implies ${b}\leq \tbigvee F$ for some finite $F\subseteq A$;
% \item  \emph{extremally disconnected} if $a^*\vee a^{**}=1$ for every$a\in L$; and
 
 \item  \emph{zero-dimensional} if each element of $L$ is a join of complemented elements.
 \end{enumerate}

\subsection{Sublocales}

	A \emph{sublocale set} (a \emph{sublocale}, for short) $S$ of a frame $L$ is a subset $S \subseteq L$ such
	that
	\begin{itemize}
		\item[(S1)] for every $A \subseteq S$, $\bigwedge A$ is in $S$, and
		\item[(S2)] for
		every $s \in S$ and every $x \in L$, $x \to s$ is in $S$.
	\end{itemize}
	The system of sublocale sets ordered by inclusion, $\mathsf{S}(L)$, is a coframe (meet coincides with the intersection and the join  is given by $\tbigvee S_i=\{\tbigwedge M\mid M\subseteq \tbigcup S_i\};$ the top is $L$ and the bottom is the set $\mathsf{O}=\{1\}$).

\subsubsection{Closed and open sublocales} To each $a\in L$ one associates two
sublocales,
\[
\mathfrak{c}(a)=\newuparrow a
\qquad\text{and}\qquad
\mathfrak{o}(a)=\{a\rightarrow b\mid b\in L\},
\]
called \emph{closed} and \emph{open}\index{open sublocale $\sim$} respectively. The
assignment $a\mapsto\mathfrak{c}(a)$ is a coframe embedding
$L^{\mathrm{op}}\hookrightarrow\mathsf{S}(L)$; it identifies $L^{\mathrm{op}}$ with the
subcoframe $\mathfrak{c}(L)\subseteq\mathsf{S}(L)$ of all closed sublocales. Writing
$\mathfrak{o}(L)$ for the subcoframe of $\mathsf{S}(L)$ generated by the
$\mathfrak{o}(a)$, the assignment $a\mapsto\mathfrak{o}(a)$ is instead a poset embedding
$L\rightarrow\mathfrak{o}(L)$.

The following identities hold in $\mathsf{S}(L)$:
\begin{enumerate}[(1)]
\item $\mathfrak{o}(a)$ and $\mathfrak{c}(a)$ are mutually complementary.
\item $\mathfrak{c}(0)=1$, $\mathfrak{c}(1)=0$,
$\mathfrak{c}(a)\vee\mathfrak{c}(b)=\mathfrak{c}(a\wedge b)$ and
$\tbigwedge_{i\in I}\mathfrak{c}(a_i)=\mathfrak{c}\!\left(\tbigvee_{i\in I}a_i\right)$.
\item $\mathfrak{o}(0)=0$, $\mathfrak{o}(1)=1$,
$\mathfrak{o}(a)\wedge\mathfrak{o}(b)=\mathfrak{o}(a\wedge b)$ and
$\tbigvee_{i\in I}\mathfrak{o}(a_i)=\mathfrak{o}\!\left(\tbigvee_{i\in I}a_i\right)$.
\item Every $S\in\mathsf{S}(L)$ admits a representation
$S=\tbigwedge_i(\mathfrak{o}(a_i)\vee\mathfrak{c}(b_i))$.
\end{enumerate}

\subsubsection{Subspaces and induced sublocales} A subspace $A$ of a $T_0$-space $X$
gives rise to a sublocale of $\mathcal{O}X$, formed by those $U\in\mathcal{O}X$ that
satisfy
\[
U=\tbigcup\{V\in\mathcal{O}X\mid V\cap A=U\cap A\};
\]
we write $\widetilde{A}$ for it. For an open set $U$ this recovers the open and closed
sublocales,
\[
\widetilde{U}=\mathfrak{o}(U)\quad\text{and}\quad
\widetilde{X\setminus U}=\mathfrak{c}(U),
\]
and, more generally, the construction turns unions into joins:
\[
\widetilde{\tbigcup_{i\in I}A_i}=\tbigvee_{i\in I}^{\mathsf{S}(L)}\widetilde{A_i}.
\]
Thus a join of induced sublocales is again induced. The same is not true of meets:
in general $\widetilde{A\cap B}\neq\widetilde{A}\wedge\widetilde{B}$.

\subsubsection{Booleanization}\label{booleanization} Each frame $L$ carries an
associated complete Boolean algebra $\mathsf{B}(L)$, whose elements are the regular
elements $a=a^{\ast\ast}$, together with a frame homomorphism
\[
\beta_L\colon L\to\mathsf{B}(L),\qquad a\mapsto a^{\ast\ast}.
\]
This $\mathsf{B}(L)$ is a sublocale of $L$. Call a sublocale $S$ \emph{dense} when $0\in S$. Isbell's density theorem asserts that
among the dense sublocales of $L$ there is a least one, and that it is precisely
$\mathsf{B}(L)$.

\subsubsection{Boolean sublocales} For $a\in L$, the set
\[
\mathfrak{b}(a)=\{x\to a\mid x\in L\}
\]
is the least sublocale to which $a$ belongs, and a sublocale $S$ is Boolean exactly
when $S=\mathfrak{b}(a)$ for some $a$. Every sublocale decomposes as a join (in fact a
union) of such Boolean pieces:
\[
S=\tbigvee_{a\in S}\mathfrak{b}(a)=\tbigcup_{a\in S}\mathfrak{b}(a).
\]

\subsection{Points and the largest spatial sublocale of $L$}

An element  $p\neq 1$ is \emph{meet-prime} (briefly, \emph{prime}) if
\[
a\wedge b\le p\ \Longrightarrow\ a\le p\ \text{ or }\ b\le p
\qquad(a,b\in L).
\]
The primes of $L$ are exactly the \emph{points} of the locale: they correspond bijectively to the frame homomorphisms $\xi\colon L\to\{0,1\}$, a point $\xi$ being
recovered as $p_\xi=\tbigvee\{a\in L\mid \xi(a)=0\}$. We write $\mathsf{Prm}(L)$ for the set of primes of $L$.

For each prime $p$ we have that
$
\{p,1\},
$
is a sublocale, precisely the smallest sublocale containing $p$, namely, $\mathfrak{b}(p)$ (indeed, for $a\not\le p$ one has
$a\rightarrow p=p$, while $a\le p$ gives $a\rightarrow p=1$). As a frame $\mathsf{b}(p)$
is the two-element Boolean algebra, so it represents a single point; accordingly we
call it the \emph{one-point sublocale} determined by $p$. The one-point sublocales are
exactly the atoms of $\mathsf{S}(L)$.

Joining all of them yields the sublocale
\[
\mathsf{Sp}(L)=\tbigvee_{p\in\mathsf{Prm}(L)}\mathfrak{b}(p),
\]
the \emph{largest spatial sublocale} of $L$: every spatial sublocale is a join of one-point sublocales, hence is contained in $\mathsf{Sp}(L)$, and $\mathsf{Sp}(L)$ is itself spatial. This is precisely isomorphic to the image of the spatial reflection, namely, $\mathsf{Sp}(L)\cong\Omega\Sigma L$. In particular, $L$ is spatial if and only if $\mathsf{Sp}(L)=L$, that is, if and only if each element of $L$ is a meet of primes. Obviously $\mathsf{Sp}(L)$ is the largest sublocale of $L$ containing all of its points.

\section{Mereological sums and exact joins}\label{sums.exactjoins}

In this section we study the relation between mereological sums and exact joins and show that both notions coincide in lattices $P$ such that $P^+$ is separative. This analysis will be useful later to better understand the Boolean completion of $\mathsf{S}(L)$ as its sum-completion or exact join-completion.

\begin{proposition}

Let $A$ be a nonempty subset of a lattice $P$. If $b=\tbigvee A$ is exact, then $b\msum A$ in $P^+$.
\end{proposition}

\begin{proof}If $b=0$, it is obvious. Otherwise,
let $0\neq c\leq b$ in $P$. Then, by Theorem \ref{exactdist},
\[
0\neq c=c\wedge b=\tbigvee_{a\in A} (c\wedge a),
\]
thus there exists $a\in A$ such that $c\wedge a\neq 0$, that is $c\circ a$ in $P^+$. Therefore $b\msum A$ in $P^+$.
\end{proof}

The following result follows directly from Lemma 4 and Lemma 9 in \cite{GG13}, but we provide here a direct proof for completeness.

\begin{proposition}
Let $P$ be a separative poset, $a,b\in P$ and $B\subseteq P$ such that $a\msum B$ and $b\msum B$. Then $a=b$.
\end{proposition}

\begin{proof}
Assume that $a\not\leq b$. Since $P$ is separative, there exists $c\in P$ such that $c\leq a$ and $c\wr b$. Then $c\wr x$  for any $x\in B$ as $B\leq b$, but we also have that $c\circ a$ and $a\msum B$, a contradiction. Therefore $a\leq b$. Dually, we also have that $b\leq a$.  
\end{proof}

\begin{proposition}
Let $P$ be a lattice such that $P^+$ is separative. If
$a\msum B$ in $P^+$ then $a$ is the exact join of $B$.
\end{proposition} 
\begin{proof}
Let $c<d$  in $P$ such that $B\subseteq c\downarrow d$. If $a\wedge d\not\leq c$, as $P^+$ is separative, there exists $x\in P^+$ such that  $ x\leq a\wedge d$ and $x\wr c$, that is, $x\wedge c=0$. For any $b\in B$,
\[
x\wedge b=x\wedge b\wedge d 
\]
since  $x\leq d$ and therefore 
\[
x\wedge b\leq x\wedge c=0
\]
as $b\wedge d\leq c$ since $B\subseteq c\downarrow d$. In consequence we have that $x\circ a$ while $x\wr b$ for all $b\in B$, a contradiction since $a$ is a sum of $B$. We conclude that $a\in c\downarrow d$. Thus $a$ is the exact join of $B$.
\end{proof}

\begin{corollary}\label{exactjoins.equals.sums}
Let $P$ be a lattice such that $P^+$ is separative. Then sums in $P^+$ and exact joins in $P$ distinct from 0 coincide.
\end{corollary}

\section{The Boolean completion of $\mathsf{S}(L)$}\label{section.BooleanCompletion}

Consider the topology $\mathsf{R}(L)$ on $L$ generated by $\mathsf{S}(L)$  as a basis. In $\mathsf{R}(L)$ every element $x$ of $L$ has a smallest neighborhood, namely, $\mathfrak{b}(x)$. This is equivalent to say that $\mathsf{R}(L)$ is an Alexandrov topology, a topology in which arbitrary intersections of open susbsets are open. Accordingly, we have the following:

\begin{lemma}
Let $S\subseteq L$. Then $S\in \mathsf{R}(L)$ if and only if
\[
S=\tbigcup_{a\in S}\mathfrak{b}(a).
\]
\end{lemma}

\begin{proof}
Obvious.
\end{proof}

\begin{definition}
We will say that two elements $x$ and $y$ of a frame $L$ are \emph{apart} and denote it by $x\perp y$ if there are sublocales separating them, that  is, if there exist $S,T\in\mathrm{S}(L)$ such that $x\in S$, $y\in T$ and $S\cap T=\{1\}=\mathsf{O}$ (namely, if $x$ and $y$ can be separated by open subsets in $L\setminus\{1\}$ endowed with the subspace topology induced by $\mathsf{R}(L)$). Also, for subsets $U\subseteq L$, we will write $U\perp y$ to mean that $x\perp y$ holds for all $x\in U$.
\end{definition}

\begin{lemma} Let $x, y\in L$ and $S\in\mathsf{S}(L)$. 
	\begin{enumerate}[\rm(1)]
		\item $x\perp y$ if and only if $\mathfrak{b}(x)\wedge\mathfrak{b}(y)=\mathsf{O}$.
		\item $S\perp x$ if and only if $S\wedge \mathfrak{b}(x)=\mathsf{O}$.
	\end{enumerate}
\end{lemma}

\begin{proof}
		\begin{enumerate}[\rm(1)]
		\item Obvious.
		\item If $S\perp x$, by (1), we have that $\mathfrak{b}(y)\wedge\mathfrak{b}(x)=\mathsf{O}$ for all $y\in S$. As $S=\tbigcup_{y\in S}\mathfrak{b}(y)$, we have that 
		\[ S\wedge \mathfrak{b}(x)=\left(\tbigcup_{y\in S}\mathfrak{b}(y)\right)\cap \mathfrak{b}(x)
		=\tbigcup_{y\in S}\left(\mathfrak{b}(y)\cap\mathfrak{b}(x)\right)=\mathsf{O}.
		\]
		The other implication is obvious.\qedhere
	\end{enumerate}
\end{proof}

 Note that, in $\mathsf{R}(L)$, the empty set $\varnothing$ is covered by the open subset $\{1\}$. Thus, in particular, 
 \[
 \mathsf{R}(L)^+=\mathsf{R}(L)\setminus\{{\varnothing}\}
 \]
 is still a frame.

\begin{remark}
Throughout the remainder of the paper, unless explicitly indicated otherwise, pseudocomplements and complements of elements of $\mathsf{R}(L)^+$ will always be understood as those computed in $\mathsf{R}(L)^+$. Accordingly, if $S\in\mathsf{R}(L)^+$, the notation $S^\ast$ and $S^c$ refers to the pseudocomplement and the complement of $S$ in $\mathsf{R}(L)^+$, respectively, regardless of whether S is a sublocale set.
\end{remark}

%\begin{definition}
%Let $\mathsf{R}(L)^+$ denote the collection of non-empty unions of sublocales sets of a frame $L$.
%\end{definition}
%
%\begin{lemma}
%$\mathsf{R}(L)^+$ ordered by inclusion is a frame.
%\end{lemma}
%
%\begin{proof}
%Simply note that as $\mathsf{S}(L)$ is closed under intersections, $\mathsf{R}(L)^+\cup\{\varnothing\}$ is a topology on $L$, therefore a frame ordered by inclusion. This frame has a single element $\{1\}$ above the bottom element and below all other elements. Thus  just dropping the bottom element $\{\varnothing\}$ still yields a frame.
%\end{proof}

\begin{proposition}\label{regular.elements}
Regular elements of $\mathsf{R}(L)^+$ are precisely nonempty subsets $S\subseteq L$  satisfying that $x\in S$ whenever $x\perp  y$  for  all $y\in L$ with $S\perp y$. 
	
\end{proposition}

\begin{proof}
	First note that any nonempty subset $S\subseteq L$ satisfying the condition is an element of $\mathsf{R}(L)^+$.  Indeed, let $a\in S$, $b\in \mathfrak{b}(a)$ and $y\in L$ such that $S\perp y$. Then $a\perp y$. Consequently
\[
\mathsf{O}=\mathfrak{b}(a)\cap \mathfrak{b}(y)\supseteq \mathfrak{b}(b)\cap\mathfrak{b}(y).
\]
Therefore $b\perp y$ and, by hypothesis, $b\in S$. Thus $\mathfrak{b}(a)\subseteq S$. Consequently we have
\[
S=\tbigcup_{a\in S}\mathfrak{b}(a),
\]
that is, $S\in \mathsf{R}(L)^+$.
	
	In order to check that the condition characterizes regular elements,  let $S\in \mathsf{R}(L)^+$. Then one has
	\[
	\begin{aligned}
	S^\ast&=\tbigcup\{T\in\mathsf{S}(L)\mid S\cap T=\mathsf{O}\}\\
		&=\tbigcup\{\mathfrak{b}(y)\mid S\perp y\}\\
		&=\{y\mid S\perp y\}\\
	\end{aligned}
	\]
	Consequently,
	\[
	\begin{aligned}
		S^{\ast\ast}&=\{z\mid S^\ast\perp z\}\\
			&=\{z\mid y\perp z \text{ whenever } S\perp y\}.
	\end{aligned}
	\]
	If $S$ satisfies the condition, then clearly $S=S^{\ast\ast}$.
	For a general $S$, as $S^{\ast\ast\ast\ast}=S^{\ast\ast}$, we have that
	\[
	S^{\ast\ast}=\{z\mid y\perp z \text{ whenever }S^{\ast\ast}\perp y\},
	\]
	thus $S^{\ast\ast}$ satisfies the condition.
\end{proof}

\begin{remark}
	Note that if we drop the top element $1$ from $L$ and consider the subspace topology induced by $\mathsf{R}(L)$, we obtain a lattice of open subsets isomorphic to $\mathsf{R}(L)^+$. Moreover, we can easily modify Proposition \ref{regular.elements} to obtain a characterization of regular open subsets of arbitrary Alexandrov-discrete spaces.
\end{remark}

In \cite{Wallman38}, Wallman applied lattice-theoretic methods to prove topological results, requiring a point-free replacement for the $T_1$ separation axiom. His notion of \emph{disjunctivity} worked remarkably well. It is usually presented as 
\[
a\not \leq b \implies \exists c, a\wedge c\neq 0=b\wedge c.
\]
Note that a lattice $L$ is disjunctive if and  only if $L^+$ is separative.
Its dual form is called \emph{conjunctivity} or \emph{subfitness} \cite{I72} and 
$\mathsf{S}(L)^{\mathrm{op}}$ is well known to be subfit as it is zero-dimensional \cite{PPT19}. So, in other words, $\mathsf{S}(L)$ is well known to be disjunctive and we will use this fact in the following proof.

\begin{proposition}\label{sublocales.are.parts}
	Sublocales are regular elements of $\mathsf{R}(L)^+$.
\end{proposition}

\begin{proof}
	For any sublocale $S$ of $L$, we always have that $S\subseteq S^{\ast\ast}$, so we only have to check that $S^{\ast\ast}\subseteq S$. Assume that $S^{\ast\ast}\not\subseteq S$. Then, as $S^{\ast\ast}$ is a union of sublocale sets, there exists $T\in\mathsf{S}(L)$ such that $T\subseteq S^{\ast\ast}$ and $T\not\subseteq S$. By disjunctivity, there exists $R\in \mathsf{S}(L)$ such that
	\[
	S\cap R=\mathsf{O}\neq T\cap R.
	\]
	As $T\subseteq S^{\ast\ast}$, we have $S^{\ast\ast}\cap R\neq \mathsf{O}$.	Let $1\neq x\in S^{\ast\ast}\cap R$. Then
	\[
	S\cap\mathfrak{b}(x)=\mathsf{O}\neq S^{\ast\ast}\cap \mathfrak{b}(x).
	\]
	Accordingly we have that $x\in S^{\ast}$ while $x\not\in S^{\ast\ast\ast}$, a contradiction, since $S^{\ast}=S^{\ast\ast\ast}$.
\end{proof}

\begin{theorem}
The inclusion map $\mathsf{S}(L)^+\to \mathsf{B}(\mathsf{R}(L)^+)$, where $\mathsf{B}(\mathsf{R}(L)^+)$ is the Booleanization of $\mathsf{R}(L)^+$, constitutes the Boolean completion of $\mathsf{S}(L)^+$.
\end{theorem}

\begin{proof}
	It is obvious that $\mathsf{S}(L)\setminus\{0\}$ is dense in $\mathsf{B}(\mathsf{R}(L)^+)$.
\end{proof}

It was explained in \cite{Mozo2024} that the Boolean completion is, up to isomorphism, the only sum-completion once strong supplementation is assumed. Now, taking into account Corollary \ref{exactjoins.equals.sums}, we can reinterpret the Boolean completion $\mathsf{S}(L)^+\to \mathsf{B}(\mathsf{R}(L)^+)$ in the theorem above as the \emph{exact-join-completion} of $\mathsf{S}(L)^+$ in the following sense: this is the only order-embedding of $\mathsf{S}(L)^+$ it into a lattice $M$ such that $M^+$ is separative, all its non-empty subsets have a exact-join and  all elements $a\in M^+$ can be obtained as the exact-join of the image under the order-embedding of a subset of $\mathsf{S}(L)^+$.

\section{Parts of a frame $L$} \label{section.Parts}

\begin{definition}
We call the elements of $\mathsf{B}(\mathsf{R}(L)^+)$ \emph{parts} of $L$ and denote the Boolean algebra of all parts of $L$, namely, $\mathsf{B}(\mathsf{R}(L)^+)$, by $\mathsf{P}(L)$.
\end{definition}

\begin{remark}
As  infima in $\mathsf{R}(L)^+$ are given by intersections, the inclusion map $j_L$ of $\mathsf{S}(L)$ into $\mathsf{R}(L)^+$ preserves arbitrary infima. Since $\beta_{\mathsf{R}(L)^+}$ also preserves finite infima, we have that $e_L=\beta_{\mathsf{R}(L)^+}\circ j_L$ also does so. Now we have the situation depicted in the following diagram
\[
\xymatrix{\mathsf{S}(L)\ar[r]^{j_l}\ar@<1ex>[dr]^{e_L}&\mathsf{R}(L)^+\ar@{>>}[d]^{\beta_{\mathsf{R}(L)^+}}\\
	&\mathsf{B}(\mathsf{R}(L)^+)=\mathsf{P}(L)\ar@<1ex>[ul]^{\tau_L}
}
\]
with $\tau_L$ defined as
\[
\tau_L(P)=\tbigvee^{\mathsf{S}(L)}\{S\in\mathsf{S}(L)\mid S\leq P\}
\]
for each $P\in\mathsf{P}(L)$.

\end{remark}

\begin{lemma}\label{tau.infima}

Let $P\in\mathsf{P}(L)$. Then
\[
\tau_L(P)=\tbigwedge^{\mathsf{S}(L)}\{S\in\mathsf{S}(L)\mid S\geq P\}%\left(=\tbigcap^{\mathsf{S}(L)}\{S\in\mathsf{S}(L)\mid S\supseteq P\}\right)
.
\]

\end{lemma}

\begin{proof}
It is obvious that
\[
\tau_L(P)\leq \tbigwedge^{\mathsf{S}(L)}\{S\in\mathsf{S}(L)\mid S\geq P\}.
\]
 Also, $P\leq \tau_L(P)$ since  $\mathfrak{b}(a)\leq P$ for each $a\in P$. Since $\tau_L(P)$ is a sublocale, we conclude that 
\[
\tau_L(P)\geq \tbigwedge^{\mathsf{S}(L)}\{S\in\mathsf{S}(L)\mid S\geq P\}
\]
for any $P\in\mathsf{P}(L)$.
\end{proof}

\begin{lemma}$\tau_L$ is the left adjoint to $ e_L$, 
\[
\tau_L\dashv e_L,
\]
i.e., for any $P\in\mathsf{P}(L)$ and $S\in\mathsf{S}(L)$, we have
\[
P\leq e_L(S) \iff  \tau_L(P)\leq S.
\]
\end{lemma}

\begin{proof}
Note that, by Proposition \ref{sublocales.are.parts}, $e_L$ is simply the inclusion map. Therefore we have $e_L(S)=S$. We also have that $\tau_L(S)$=S. As $\tau_L$ is clearly monotone and, by Lemma \ref{tau.infima}, $\tau_L(P)\leq S$  implies that $P\leq S$, we have that 
\[
P\leq e_L(S)=S
\iff \tau_L(P)\leq\tau_L(S)=S.\qedhere
\]
\end{proof}

\begin{corollary}
$\tau_L$ preserves arbitrary suprema and $e_L$ preserves arbitrary infirma.
\end{corollary}

\begin{proposition}\label{complemented.sublocales}
Let $S$ and $T$ be sublocales of $L$ complements to each other in $\mathsf{S}(L)$. Then they are also complements to each other as parts, that is, as elements of $\mathsf{P}(L)$.
\end{proposition}

\begin{proof}
As $e_L$ preserves finite infima, we have that 
\[
0=S\wedge_{\mathsf{S}(L)}T=e_L(S)\wedge_{\mathsf{P}(L)}e_L(T)=S\wedge_{\mathsf{P}(L)}T.
\]
Now, let $ R\in\mathsf{S}(L)$ such that
\[
R\wedge_{\mathsf{P}(L)}(S\vee_{\mathsf{P}(L)}T)=0.
\]
As $S\cup T\subseteq S\vee_{\mathsf{P}(L)}T$ and $\wedge_{\mathsf{P}(L)}=\cap$, we have

\[
R\cap(S\cup T)=(R\cap S)\cup(R\cap T)=0.
\]
 As $S$ and $T$ are complement to each other in $\mathsf{S}(L)$, we have that $R\leq T$ and $R\leq S$. Therefore $R\leq S\cap T=0$. As $\mathsf{P}(L)$ is Boolean and its elements are unions of sublocale sets, we conclude that $S\vee_{\mathsf{P}(L)}T=1$.
\end{proof}

\begin{proposition}
	The map $\iota_L\colon L\to\mathsf{P}(L)$ given by $a\mapsto \mathfrak{o}(a)$ is a frame embedding.
\end{proposition}

\begin{proof}As  $\iota_L$ is clearly injective, we only need to check that it preserves finite meets and arbitrary joins. For this end, let $a, b\in L$. As $\mathfrak{o}(a)\cap\mathfrak{o}(b)
	=\mathfrak{o}(a\wedge b)$ is in $\mathsf{P}(L)$, we have that 
	\[
	\mathfrak{o}(a)\wedge_{\mathsf{P}(L)}\mathfrak{o}(b)=\mathfrak{o}(a\wedge b).
	\]
	Accordingly, $\iota_L$ preserves finite meets.
	
	Now let $\{a_i\}_{i\in I}\subseteq L$ and note that
	\[
	\begin{aligned}
		\left(\tbigcup_{i\in I}\mathfrak{o}(a_i)\right)^\ast
		&=\tbigcup\left\{T\in\mathsf{S}(L)\mid T\cap\tbigcup_{i\in I} \mathfrak{o}(a_i)=\mathsf{O}\right\}\\
		&=\tbigcup\left\{T\in\mathsf{S}(L)\mid \forall i\in I, T\cap \mathfrak{o}(a_i)=\mathsf{O}\right\}\\
		&=\tbigcup\left\{T\in\mathsf{S}(L)\mid \forall i\in I, T\subseteq \mathfrak{c}(a_i)\right\}\\
		&=\tbigcup\left\{T\in\mathsf{S}(L)\mid  T\subseteq \mathfrak{c}\left(\tbigvee_{i\in I}a_i\right)\right\}\\
		&=\mathfrak{c}\left(\tbigvee_{i\in I}a_i\right).
	\end{aligned}	
	\]
	As, by Proposition \ref{complemented.sublocales},  we have that
	\[
	\mathfrak{c}\left(\tbigvee_{i\in I}a_i\right)^\ast=\mathfrak{o}\left(\tbigvee_{i\in I}a_i\right),
	\]
	we conclude that 
	\[
	\left(\tbigcup_{i\in I}\mathfrak{o}(a_i)\right)^{\ast\ast}=\mathfrak{o}\left(\tbigvee_{i\in I}a_i\right).
	\]
	Therefore
	\[
	\left(\tbigvee_{i\in I}^{\mathsf{B}(\mathsf{R}(L)^+)}\mathfrak{o}(a_i)\right)^{\ast\ast} = \mathfrak{o}\left(\tbigvee_{i\in I}a_i\right)
	\]
	and we conclude that $\iota_L$ preserves arbitrary joins. 
\end{proof}

Therefore we have an isomorphic copy of $L$ in $\mathsf{P}(L)$,
\[
\mathfrak{o}L=\{\mathfrak{o}(a)\mid a \in L\}.
\]

\section{Some interesting parts}\label{section.interesting}

\subsection{The pointless part of a frame}

A set $\{1,m\}$ with $1\neq m\in L$ is a sublocale of a frame $L$ if and only if $m$ is prime in $L$. These are called \emph{one-point sublocales} \cite{PP2008}. Obviously, $\mathfrak{b}(m)=\{1,m\}$ in this case. We will call
\[
\mathrm{Points}(L)=\tbigvee^{\mathsf{P}(L)}\{\mathfrak{b}(m)\mid m\text{ is prime in } L \}
\]
 the \emph{part of merely points} of $L$ and 
and  its complement in $\mathsf{P}(L)$,
\[
\mathrm{Pointless}(L)=\mathrm{Points}(L)^c,
\]
 the \emph{pointless part}  of $L$.

\begin{lemma}\label{meet}
Let $a,m\in L$ such that $m$ is prime. If $m\in \mathfrak{b}(a)$, there exists $b\in L$ such that $a\leq b$  and $a$ is prime in $\downarrow\!b$.
\end{lemma}

\begin{proof}
	If $m\in\mathfrak{b}(a)$ there exists $b\in L$ such that $b\to a=m$. Let 
	\[
	b^\prime=\tbigvee\{x\mid x\to a\geq m\}.
	\]
	As $a\to a=1\geq m$, we have that $b^\prime\geq a$. We also have that
	\[
	\begin{aligned}
	b^\prime\to a&=
		\left(\tbigvee\left\{x\mid x\to a\geq  m\right\}\right)\to a\\
		&=\tbigwedge\left\{x\to a\mid x\to a\geq  m\right\}=m.
	\end{aligned}
	\]
	Now let $y, z\leq b^\prime$ such that $y\wedge z\leq a$. Since $a\leq b^\prime\to a=m$, without loss of generality we can assume that $y\leq m$. Then
	\[
	\begin{aligned}
	y&\leq b^\prime \wedge m\\
	&= b^\prime \wedge (b^\prime \to a)\leq a.
	\end{aligned}
	\]
	We conclude that $a$ is prime in $\downarrow\! b^\prime$.
\end{proof}
	
\begin{lemma}\label{neces}
Let $a,b\in L$ such that $a\leq b$. Then $a$ is prime in $\downarrow\! b$ if and only if $b\to a$ is prime in $L$.
\end{lemma}

\begin{proof}
	Let $x, y\in L$ such that $x\wedge y \leq b\to a$. Then
	\[
	(x\wedge b)\wedge (y\wedge b)=x\wedge y\wedge b\leq a.
	\]
	Then, if $a$ is prime in $\downarrow\! b$, without loss of generality, we can assume that $x\wedge b\leq a$. Thus $x\leq b\to a$. We conclude that $b\to a$ is prime.
	
	Now let $x,y\in \downarrow\! b$ such that $x\wedge y\leq a$. As $a\leq b\to a$, we have that $x\wedge y\leq b\to a$. If $b\to a$ is prime, then either $x\leq b\to a$ or $y\leq b\to a$. Say that the former holds. Then $x=x\wedge b\leq a$. Accordingly, $a$ is prime in $\downarrow\! b$.
\end{proof}

\begin{definition}
	We will say that $a\in L$ is \emph{relatively prime} if there exists $b\geq a$ in $L$ such that $a$ is prime in $\downarrow\! b$. 
\end{definition}

\begin{lemma}
	Let $a\in L$. Then $\mathfrak{b}(a)$ contains a relatively prime element if and only if $a$ is relatively prime. 
\end{lemma}

\begin{proof}
If $a$ is relatively prime then $\mathfrak{b}(a)$ obviously contains a relatively prime element, $a$ itself. If $x\in\mathfrak{b}(a)$ is relatively prime, then there is a prime $m\in L$ such that \[
m\in\mathfrak{b}(x)\subseteq\mathfrak{b}(a).
\]
 By Lemma \ref{meet}, we have that $a$ is relatively prime.
\end{proof}

\begin{theorem}\label{pointless}
	Let $a\in L$. Then $a\in\mathrm{Pointless}(L)$ if and only $a$ is not relatively prime.
\end{theorem}

\begin{proof}
Note that $a\in\mathrm{Pointless}(L)$ if and only if, for all prime elements $m$,
$\mathfrak{b}(a)\cap\mathfrak{b}(m)=\mathsf{O}$, that is, $m\not\in\mathfrak{b}(a)$. Then this follows easily from Lemma \ref{neces}.
\end{proof}

Note that a related notion appears in work of Ball  (\cite{B24}), who uses the term \emph{pointless part} in the context of sublocales. In that setting, the construction is formulated purely in terms of sublocales of a frame. In the 
present paper, however, the notion of a \emph{part} of a frame is introduced as a specific structure, and the pointless part of a frame is defined within this more general framework. Thus our description of the pointless part applies in a broader context. Nevertheless, when one restricts attention to $T_1$-frames, namely, frames in which prime elements are precisely maximal elements, our pointless part coincides with Ball's pointless sublocale. The following example shows that our result is a proper generalizations of Ball's, as, in general, $\mathrm{Pointless}(L)$ is not a sublocale if the underlying frame $L$ is no $T_1$.
	
\begin{example} %$\mathrm{Pointless}(L)$ is not a sublocale in general. As a counterexample, 
Let $L=\mathcal{N}_0\cup\{\varnothing\}$ where $\mathcal{N}_0$ is the filter of open neighborhoods of 0 of the real line $\mathbb{R}$ endowed with the Euclidean topology. Clearly $L$ is a frame. It is easy to check  that relatively prime elements in this frame are either the empty set or open neighborhoods of $0$ that have a successor under the inclusion order, that is, a relatively prime element $V$ is either the empty set or of the form $V=U\setminus\{x\}$ where  $U$ is an open neighborhood  of 0 and $0\neq x\in U$. Then we have that $\bigcap_{p<0<p}(p,q)=\varnothing$ is a prime element of $L$ while $(p,q)$ with $p<0<q$ are not relatively prime elements. Therefore, $\mathrm{Pointless}(\mathcal{N}_0)$ is not closed under infima in $L$ and, consequently, neither a sublocale.
\end{example}

\subsection{The part of merely points}

\begin{lemma}\label{points.all.points}
	$a\in\mathrm{Points}(L)$ if and only if, for  all $1\neq x\in\mathfrak{b}(a)$, $\mathfrak{b}(x)$ contains a prime element of $L$.
\end{lemma}

\begin{proof}
	By Theorem \ref{pointless}, we have that $a\in\mathrm{Points}(L)$ if and only if $\mathfrak{b}(a)$ has not element which is not relatively prime, that is, if an only if, $\mathfrak{b}(a)$ contains only relatively prime elements and 1. By Lemmata \ref{meet}  and \ref{neces}, this is equivalent to saying that, for all $1\neq x\in \mathfrak{b}(a)$, $\mathfrak{b}(x)$ contains a prime.
\end{proof}

A collection $\{x_i\}_{i\in I}\subseteq L$ is said to a be \emph{co-discrete} system if there exists a cover $C$ such that for each $c\in C$, $c\leq x_i$ for all $i$ with possibly one exception. For any co-discrete system $\{x_i\}_{i\in I}$ and any $y\in L$, we have that $y\vee\tbigwedge_{i\in I}x_i=\tbigwedge_{i\in I} y\vee x_i$ \cite{P84}.
	
\begin{theorem}\label{char.points}
	
$a\in\mathrm{Points}(L)$ if and only if $a=\tbigwedge m_i$ where $\{m_i\}_{i\in I}$ is a co-discrete system of prime elements of $L$. 
\end{theorem}

\begin{proof}
	Let $a=\tbigwedge_{i\in I}m_i$ where $\{m_i\}_{i\in I}$ is a co-discrete system of prime elements of $L$. Let $C$ be the cover that witnesses the co-discreteness of  $\{m_i\}_{i\in I}$. Note that for each $i\in I$, there exists some $c\in C$ such that $c\not\leq m_i$. Let
\[
c_i=\tbigvee\{c\in C\mid c\not\leq m_i\}.
\]
 Note also that $c_i\leq m_j$ for all $j\neq i$ in I.
	
	Now let $x\in L$. If $x\leq a$, then $x\rightarrow a=1$. Otherwise, there exists $i\in I$ such that $x\not\leq m_i$. As $m_i$ is prime, $c_i\wedge x\not \leq m_i$. Since 
\[
m_i\wedge c_i\wedge x\leq \tbigwedge_{j\in I}m_j=a,
\]
 we have that 
\[
m_i\leq (c_i\wedge x)\rightarrow a.
\]
 For any $y\in L$ such that 
\[
y\wedge c_i\wedge x\leq a=\tbigwedge_{j\in I}m_j,
\]
 we  have that 
\[
y\wedge c_i\wedge x\leq m_i.
\]
 As $m_i$ is prime and $c_i\wedge x\not\leq m_i$, we conclude that $y\leq m_i$. Therefore 
\[
m_i=(c_i\wedge x)\rightarrow a=c_i\rightarrow(x\rightarrow a),
\] that is, $m_i\in \mathfrak{b}(x\rightarrow a)$.
	
	Now let $a\in\mathrm{Points}(L)$. Let $\{m_i\}_{i\in I}$ be the collection  of all prime  elements contained in $\mathfrak{b}(a)$. Assume that $a<\tbigwedge_{i\in I} m_i$. Then we have that 
\[
\mathfrak{b}(a)\not\leq\mathfrak{b}\left(\tbigwedge_{i\in I}m_i\right).
\]
 By disjunctivity, there exists $c\in L$ such that 
\[
\mathsf{O}\neq \mathfrak{b}(c)\subseteq \mathfrak{b}(a)
\]
 and 
\[	
\mathfrak{b}(c)\cap\mathfrak{b}\left(\tbigwedge_{i\in I}m_i\right)=\mathsf{O}.
\]
 We conclude that $\mathfrak{b}(c)$ does not contain any prime elements, a contradiction by Lemma \ref{points.all.points}. Thus, $a=\tbigwedge_{i\in I}m_i$.
\end{proof}

\begin{example}
Consider the overlapping interval topology on the closed interval $[-1,1]$ of the real line where open subsets are of the form $[-1,b), (a,b)$ and $(a, 1]$ with $-1\leq a<0<b\leq 1$ together with $[-1,1]$ and the empty set. Prime elements of the frame $L$ of open subsets of this space are the empty set $\varnothing$ and subsets of the form $[-1, b)$ and $(a,1]$ with $-1\leq a<0<b\leq 1$. Note also that $(a,b)$ is prime in $\downarrow\![-1, b)$ and $\downarrow\!(a,1]$ for any $a$ and $b$ such that $-1\leq a<0<b\leq 1$. Therefore, $\mathrm{Pointless}(L)=\mathsf{O}$. As $\mathrm{Pointless}(L)$ and $\mathrm{Points}(L)$ are complements to each other, we conclude that $\mathrm{Points}(L)=L$. Actually, we can also see this directly. Note that, for any $a$ and $b$ such that $-1\leq a<0<b\leq 1$, $\{[-1,b),(a,1]\}$ is a cover that  witnesses the co-discreteness of its own self. Therefore, $[-1,b)\cap(a,1]=(a,b)\in \mathrm{Points}(L)$.
\end{example}

\subsection{The topological glue}
Theorem \ref{char.points} is quite suggestive. The fact that removing all points does not produce the empty space can be metaphorically explained by the presence of a kind of ``topological glue'' holding those points together. For instance, consider the sublocales $\widetilde{\mathbb{Q}}$ and $\widetilde{I}$ of the frame of real numbers induced by the subspaces of rational and irrational numbers, respectively. One may informally argue that these two sublocales are so topologically intertwined that their intersection must still contain the information that binds them --- precisely that ``topological glue'' which prevents the space from falling apart.

However, this metaphor has a limitation. If we work only with sublocales, there is no way to isolate the mere collection of points of a frame without simultaneously retaining their glue. In other words, within the lattice of sublocales, points cannot be separated from the residual topological structure that connects them. The notion of \emph{a part of merely points} fills this gap. Admittedly, even this construction retains some topological information about the points, but only of a negative kind: it records which families of points fail to cluster, rather than how they are positively tied together.

We can now introduce the part $\mathrm{Glue(L)}$ of the glue of frame $L$ as the intersection---which coincides with the infimum in the lattice of parts---of the largest spatial sublocale $\mathsf{Sp}(L)$ and the pointless part of $L$, that is,
\[
\begin{aligned}
\mathrm{Glue}(L)=\mathsf{Sp}(L)\cap& \mathrm{Pointless}(L)\\
=\{a\in L\mid& a\text{ is not relatively prime}\\ &\text{and }  a=\tbigwedge_{i\in I}^Lm_i\text{ with }m_i\text{ is prime  for all }i\in I  \}.
\end{aligned}
\]

\subsection{The complement of the largest spatial sublocale}
As
\[
\mathsf{Sp}(L)^c=\{a\in L\mid \mathfrak{b}(a)\cap \mathsf{Sp}(L)=\mathsf{O}\},
\]
we have that $a\in \mathsf{Sp}(L)^c$ if and only if  for each $x\in \mathfrak{b}(a)$,
\[
x<\tbigwedge\{m\mid x\leq m\text{ and } m\text{ is prime}\}.
\]
Note that $\mathsf{Sp}(L)^c$ is not a sublocale in general. However, by Proposition \ref{complemented.sublocales}, it is so if $\mathsf{Sp}(L)$ is complemented in $\mathsf{S}(L)$. The case of $T_D$\footnote{A frame $L$ is 	$T_D$ if each prime element is \emph{covered}, i.e., for every prime element $m\in L$ there exists $a\in L$ such that $m<x$ implies $a\leq x$.} subfit frames $L$ was studied by Arrieta, Picado and Pultr in \cite{APP24}. Further note that if $\mathsf{Sp}(L)^c$ turns out to be a sublocale set of $L$, then obviously $\mathsf{Sp}(L)$ has to be complemented also in $\mathsf{S}(L)$.

\subsection{A general picture so far}
To finish this section, the  diagram  
 \begin{center}
 \begin{tikzpicture}
 \matrix (m) [matrix of math nodes, row sep=3em,
 column sep=2em]{
 & 1=L &   \\
 \mathsf{Sp}(L)& \mathrm{Points}(L)\vee\mathsf{Sp}(L)^c  & \mathrm{Pointless}(L)   \\
\mathrm{Points}(L)  & \mathrm{Glue}(L) & \mathsf{Sp}(L)^c   \\
 & \mathsf{O}=\{1\}&    \\};
 \path[-stealth]
 (m-4-2) edge (m-3-1)
 (m-4-2) edge (m-3-2)
 (m-4-2) edge (m-3-3)
 (m-3-1) edge (m-2-1)
 (m-3-1) edge (m-2-2)
 (m-3-2) edge (m-2-1)
 (m-3-2) edge (m-2-3)
 (m-3-3) edge (m-2-2)
 (m-3-3) edge (m-2-3)
 (m-2-1) edge (m-1-2)
 (m-2-2) edge (m-1-2)
 (m-2-3) edge (m-1-2)
 (m-3-2) edge[-,line width=6pt,draw=white] (m-2-1) edge (m-2-1)
 (m-3-2) edge[-,line width=6pt,draw=white] (m-2-3) edge (m-2-3);
\end{tikzpicture}
\end{center}
 summarizes the sublattice of $\mathsf{P}(L)$ generated by the parts of a frame $L$ that we have considered so far in this section. Here each arrow represents an inclusion, which is strict in the general case.

\section*{Declaration of generative AI and AI-assisted technologies in the manuscript preparation process. }

During the preparation of this work the authors used Claude Opus 4.8 in order assist with language editing, rephrasing, and proofreading of the text. After using this tool/service, the authors reviewed and edited the content as needed and take full responsibility for the content of the published article.

\end{document}